\documentclass[11pt]{article}
\usepackage{amsfonts,amscd,amssymb}
\newtheorem{definition}{Definition}[section]
\newtheorem{lemma}[definition]{Lemma}
\newtheorem{proposition}[definition]{Proposition}

\newtheorem{remark}[definition]{Remark}

\newtheorem{theorem}[definition]{Theorem}

\def\va{\varepsilon}

\def\ra{\rightarrow}
\def\a{\alpha}
\def\b{\beta}

\def\r{\rho}
\def\cd{\cdot}

\renewcommand{\theequation}{\thesection.\arabic{equation}}

\def\rawo\lonra{\longrightarrow}

\def\ot{\otimes}

\newcommand{\eqref}[1]{(\ref{eq:#1})}

\newenvironment{proof}{{\it Proof.}}{\hfill $ \square $ \vskip 4mm}
\begin{document}
\title{A structure theorem for quasi-Hopf comodule algebras
\thanks{Research partially supported by the EC programme LIEGRITS, 
RTN 2003, 505078, and by the bilateral project ``New techniques in 
Hopf algebras and graded ring theory'' of the Flemish and Romanian 
Ministries of Research. The first author was also partially supported by  
the programme CERES of the Romanian Ministry of Education and Research, 
contract no. 4-147/2004.}}
\author
{Florin Panaite\\
Institute of Mathematics of the 
Romanian Academy\\ 
PO-Box 1-764, RO-014700 Bucharest, Romania\\
e-mail: Florin.Panaite@imar.ro\\
\and 
Freddy Van Oystaeyen\\
Department of Mathematics and Computer Science\\
University of Antwerp, Middelheimlaan 1\\
B-2020 Antwerp, Belgium\\
e-mail: Francine.Schoeters@ua.ac.be}
\date{}
\maketitle
\begin{abstract}
If $H$ is a quasi-Hopf algebra and $B$ is a right $H$-comodule algebra 
such that there exists $v:H\rightarrow B$ a morphism of right 
$H$-comodule algebras, we prove that there exists a left $H$-module 
algebra $A$ such that $B\simeq A\# H$. The main difference comparing to  
the Hopf case is that, from the multiplication of $B$, which is 
associative, we have to obtain the multiplication of $A$, which 
in general is not;  
for this we use a canonical projection $E$ arising from the fact that 
$B$ becomes a quasi-Hopf $H$-bimodule.
\end{abstract}
\section*{Introduction}
If $H$ is a Hopf algebra and $B$ is a right $H$-comodule algebra with the 
property that there exists a morphism $v:H\rightarrow B$ of right 
$H$-comodule algebras, then it is well-known that $B$ is isomorphic as a 
right $H$-comodule algebra to a smash product $A\# H$, where $A$ is 
obtained as $A=B^{co (H)}$ and its multiplication is the restriction of the 
multiplication of $B$.\\
On the other hand, if $H$ is a quasi-bialgebra and $A$ is a left $H$-module 
algebra, then the smash product $A\# H$ introduced in \cite{bpv}  
becomes a right $H$-comodule algebra and the map  
$\;j:H\rightarrow A\# H, \;j(h)=1\# h$,   
is a morphism of right $H$-comodule algebras. This raises the natural 
problem of checking whether, for a quasi-Hopf algebra $H$ and a right  
$H$-comodule algebra $B$ such that there exists a morphism 
$v:H\rightarrow B$ of right $H$-comodule algebras, 
there exists a left $H$-module  
algebra $A$ such that $B\simeq A\# H$ as right $H$-comodule algebras.  
It is likely that $A$ appears as some sort of coinvariants of $B$,  
but it is clear that its multiplication {\it cannot} be obtained as the  
restriction of the one of $B$ (since $B$ is associative while in 
general $A$ is  
{\it not}), hence we need a different approach than in the  
Hopf case.\\
We first prove that $B$ becomes an object  
in the category $_H{\cal M}_H^H$ of quasi-Hopf $H$-bimodules as  
introduced in \cite{hn}. An object $M$ in $_H{\cal M}_H^H$ is endowed with a 
projection $E:M\rightarrow M$ and a concept of {\it coinvariants},  
$M^{co (H)}=E(M)$; applying this to $B$, we obtain a projection   
$E:B\rightarrow B$ and a subspace $B^{co (H)}=E(B)$. Now we define the 
vector space $A=B^{co (H)}$, with a multiplication  
defined by $a*a'=E(aa')$, for all $a, a'\in A$. The isomorphism  
$B\simeq A\# H$ follows from the structure theorem for quasi-Hopf  
$H$-bimodules, cf. \cite{hn}, and we only have to prove that it is a  
morphism of right $H$-comodule algebras.\\
An application of our structure theorem  
is that, if we have a smash product $A\# H$ for a quasi-Hopf algebra $H$,  
it provides a method to get $A$ back from $A\# H$, as 
$A\simeq (A\# H)^{co (H)}$. 
\section{Preliminaries}\label{sec1}
\setcounter{equation}{0}
We work over a field $k$. All algebras, linear spaces 
etc. will be over $k$; unadorned $\ot $ means $\ot_k$. Following
Drinfeld \cite{d}, a quasi-bialgebra is a fourtuple $(H, \Delta , 
\va , \Phi )$, where $H$ is an associative algebra with unit, 
$\Phi$ is an invertible element in $H\ot H\ot H$, and $\Delta :\
H\ra H\ot H$ and $\va :\ H\ra k$ are algebra homomorphisms
satisfying the identities
\begin{eqnarray}
&&(id \ot \Delta )(\Delta (h))=%
\Phi (\Delta \ot id)(\Delta (h))\Phi ^{-1},\label{q1}\\[1mm]%
&&(id \ot \va )(\Delta (h))=h\ot 1, %
\mbox{${\;\;\;}$}%
(\va \ot id)(\Delta (h))=1\ot h,\label{q2}
\end{eqnarray}
for all $h\in H$, and $\Phi$ has to be a normalized $3$-cocycle,
in the sense that
\begin{eqnarray}
&&(1\ot \Phi)(id\ot \Delta \ot id) (\Phi)(\Phi \ot 1)= (id\ot id
\ot \Delta )(\Phi ) (\Delta \ot id \ot id)(\Phi
),\label{q3}\\[1mm]%
&&(id \ot \va \ot id )(\Phi )=1\ot 1\ot 1.\label{q4}
\end{eqnarray}
The identities (\ref{q2}), (\ref{q3}) and (\ref{q4}) also imply
that
\begin{equation}\label{q7}
(\va \ot id\ot id)(\Phi )= (id \ot id\ot \va )(\Phi )=1\ot 1\ot 1.
\end{equation} 
The map $\Delta $ is called the coproduct or the
comultiplication, $\va $ the counit and $\Phi $ the reassociator.
We will use the version of Sweedler's sigma notation: $\Delta (h)=h_1\ot 
h_2$, and since $\Delta$ is only quasi-coassociative we adopt the 
further convention
\begin{eqnarray*}
&&(\Delta \ot id)(\Delta (h))= h_{(1, 1)}\ot h_{(1, 2)}\ot h_2,
\;\;\; (id\ot \Delta )(\Delta (h))=
h_1\ot h_{(2, 1)}\ot h_{(2,2)}, 
\end{eqnarray*}
for all $h\in H$. We will denote the tensor components of $\Phi$
by capital letters, and those of $\Phi^{-1}$ by small letters, 
namely
\begin{eqnarray*}
&&\Phi=X^1\ot X^2\ot X^3=T^1\ot T^2\ot T^3=Y^1\ot 
Y^2\ot Y^3=\cdots\\%
&&\Phi^{-1}=x^1\ot x^2\ot x^3=
t^1\ot t^2\ot t^3=y^1\ot y^2\ot y^3=\cdots
\end{eqnarray*}
The quasi-bialgebra $H$ is called a quasi-Hopf algebra if  
there exists an 
anti-a\-u\-to\-mor\-phism $S$ of the algebra $H$ and elements $\a , \b \in
H$ such that, for all $h\in H$, we have:
\begin{eqnarray}
&&S(h_1)\a h_2=\va (h)\a \mbox{${\;\;\;}$ and ${\;\;\;}$}
h_1\b S(h_2)=\va (h)\b ,\label{q5}\\[1mm]%
&&X^1\b S(X^2)\a X^3=1 %
\mbox{${\;\;\;}$ and${\;\;\;}$}%
S(x^1)\a x^2\b S(x^3)=1.\label{q6}
\end{eqnarray}
The axioms for a quasi-Hopf algebra imply that $\va (\a )\va (\b
)=1$, so, by rescaling $\a $ and $\b $, we may assume without loss
of generality that $\va (\a )=\va (\b )=1$ and $\va \circ S=\va $.\\
If $H$ is a quasi-Hopf algebra, following \cite{hn1}, \cite{hn2} 
we may define the elements  
\begin{eqnarray}
&&\hspace*{-1cm}
p_R=p^1\ot p^2=x^1\ot x^2\beta S(x^3),\;\;
q_R=q^1\ot q^2=X^1\ot S^{-1}(\alpha X^3)X^2,\label{right}
\end{eqnarray}
satisfying the relations (for all $h\in H$):
\begin{eqnarray}
&&q^1_1p^1\ot q^1_2p^2S(q^2)=1\ot 1,\;\;\;
q^1p^1_1\ot S^{-1}(p^2)q^2p^1_2=1\ot 1,\label{unu}\\
&&\hspace*{-1cm}
\Delta (h_1)p_R[1\otimes S(h_2)]=p_R[h\otimes 1],\;\;\;
[1\otimes S^{-1}(h_2)]q_R\Delta (h_1)=[h\otimes 1]q_R.\label{trei}
\end{eqnarray}
Let us record the following easy consequence of 
(\ref{q6}) (for $q=q_R=q^1\ot q^2$): 
\begin{eqnarray}
&&q^1\beta S(q^2)=1. \label{cucu}
\end{eqnarray}
Recall from \cite{hn1} the notion of comodule algebra over a
quasi-bialgebra.
\begin{definition}
Let $H$ be a quasi-bialgebra. A unital associative algebra
$B$ is called a right $H$-comodule algebra if there 
exist an algebra morphism $\r :B\ra B\ot H$ 
and an invertible element $\Phi _{\r }\in B\ot H\ot H$ 
such that:
\renewcommand{\theequation}{\thesection.\arabic{equation}}
\begin{eqnarray}
&&\Phi _{\r }(\r \ot id)(\r (b))=(id\ot \Delta 
)(\r (b))\Phi _{\r }, 
\mbox{${\;\;\;}$$\forall $ $b\in B$,}\label{rca1}\\[1mm]%
&&(1_B\ot \Phi)(id\ot \Delta \ot id)(\Phi _{\r })(\Phi 
_{\r }\ot 1_H)= (id\ot id\ot \Delta )(\Phi _{\r })(\r \ot id\ot
id)(\Phi _{\r }),\label{rca2}\\[1mm]%
&&(id\ot \va)\circ \r =id ,\label{rca3}\\[1mm]%
&&(id\ot \va \ot id)(\Phi _{\r })=(id\ot id\ot \va )(\Phi _{\r }
)=1_B\ot 1_H.\label{rca4}
\end{eqnarray}
\end{definition}
The first example of a right $H$-comodule algebra is $H$ itself, 
with $\rho =\Delta $ and $\Phi _{\rho }=\Phi $. For a right $H$-comodule 
algebra $(B, \rho , \Phi _{\rho })$ we will denote $\rho (b)=b_{(0)}\ot 
b_{(1)}$, for all $b\in B$. If $(B', \rho ', \Phi _{\rho '})$ is another 
right $H$-comodule algebra, a morphism of right $H$-comodule algebras 
$f:B\rightarrow B'$ is an algebra map such that $\rho '\circ f=
(f\ot id)\circ \rho $ and $\Phi _{\rho '}=
(f\ot id \ot id)(\Phi _{\rho })$.\\[2mm]
Suppose that $(H, \Delta , \varepsilon , \Phi )$ is a
quasi-bialgebra. If $U,V,W$ are left $H$-modules, define 
$a_{U,V,W}:(U\otimes V)\otimes W\rightarrow U\otimes (V\otimes W)$
by 
\begin{eqnarray*}
&&a_{U,V,W}((u\otimes v)\otimes w)=\Phi \cdot (u\otimes
(v\otimes w)).
\end{eqnarray*}
The category $_H{\cal M}$ of  
left $H$-modules becomes a monoidal category (see 
\cite{k}, \cite{m} for the terminology) with tensor product 
$\otimes $ given via $\Delta $, associativity constraints
$a_{U,V,W}$, unit $k$ as a trivial 
$H$-module and the usual left and right
unit constraints.\\%
Let again $H$ be a quasi-bialgebra. We say that a $k$-vector space
$A$ is a left $H$-module algebra if it is an algebra in the
monoidal category $_H{\cal M}$, that is $A$ has a multiplication
and a usual unit $1_A$ satisfying the 
following conditions: %
\begin{eqnarray}
&&(a a^{'})a^{''}=(X^1\cd a)[(X^2\cd a^{'})(X^3\cd 
a^{''})],\label{ma1}\\[1mm]%
&&h\cd (a a^{'})=(h_1\cd a)(h_2\cd a^{'}),
\label{ma2}\\[1mm]%
&&h\cd 1_A=\va (h)1_A,\label{ma3}
\end{eqnarray}
for all $a, a^{'}, a^{''}\in A$ and $h\in H$, where $h\ot a\ra
h\cd a$ is the left $H$-module structure of $A$. Following
\cite{bpv} we define the smash product $A\# H$ as follows: as
vector space $A\# H$ is $A\ot H$ (elements $a\ot h$ will be
written $a\# h$) with multiplication
given by %
\begin{equation}\label{sm1}
(a\# h)(a^{'}\# h^{'})=%
(x^1\cd a)(x^2h_1\cd a^{'})\# x^3h_2h^{'}, %
\end{equation}
for all $a, a^{'}\in A$, $h, h^{'}\in H$. Then $A\# H$ is an 
associative algebra with unit $1_A\# 1$. Moreover, by \cite{bpv}, 
$(A\# H, \rho , \Phi _{\rho })$ becomes a right $H$-comodule algebra, 
with $\rho :A\# H\rightarrow (A\# H)\ot H$, $\rho (a\# h)=(x^1\cdot a 
\# x^2h_1)\ot x^3h_2$ and $\Phi _{\rho }=(1\# X^1)\ot X^2\ot X^3$. Also, it 
is easy to see that the map $j:H\rightarrow A\# H$, $j(h)=1\# h$, is a 
morphism of right $H$-comodule algebras. \\
If $A$, $A'$ are left $H$-module algebras, a map $f:A\rightarrow A'$ is a   
morphism of left $H$-module algebras if it is multiplicative, unital and a  
morphism of left $H$-modules.\\
If $H$ is a quasi-Hopf algebra, $B$ an associative algebra and  
$v:H\rightarrow B$ an algebra map, then, following \cite{bpv}, we can 
introduce on the vector space $B$ a left $H$-module algebra structure, 
denoted by $B^v$ in what follows, for which the multiplication, unit and  
left $H$-action are:
\begin{eqnarray}
&&b\circ b'=v(X^1)bv(S(x^1X^2)\alpha x^2X^3_1)b'v(S(x^3X^3_2)),\;\; 
\forall \;b, b'\in B, \label{mult}\\
&&1_{B^v}=v(\beta ), \label{unit}\\
&&h\triangleright _vb=v(h_1)bv(S(h_2)), \;\;\forall 
\;h\in H,\; b\in B.\label{action}
\end{eqnarray}
If $H$ is a quasi-Hopf algebra and $A$ is a left $H$-module algebra, define 
the map 
\begin{eqnarray}
&&i_0:A\rightarrow A\# H, \;\;i_0(a)=p^1\cdot a\# p^2, \;\;
\forall \; a\in A, 
\end{eqnarray}
where $p=p_R=p^1\ot p^2$ is given by (\ref{right}).  
Then, by \cite{bpv}, $i_0$ becomes a   
morphism of left $H$-module algebras from $A$ to $(A\# H)^j$.
\section{The structure theorem}\label{sec2}
\setcounter{equation}{0}
We start with a lemma which is of independent interest.
\begin{lemma}\label{l1}
Let $H$ be a quasi-bialgebra and $A$ a left $H$-module with 
a multiplication. Define a multiplication on $A\ot H$ by  
\begin{eqnarray}
&&(a\ot h)(a'\ot h')=(x^1\cdot a)(x^2h_1\cdot a')\ot x^3h_2h', \label{lulu}
\end{eqnarray}
for all $a, a'\in A$ and $h, h'\in H$, and assume that this 
multiplication is associative. Then:\\
(i) The multiplication of $A$ satisfies the condition 
\begin{eqnarray*}
&&(ab)c=(X^1\cdot a)((X^2\cdot b)(X^3\cdot c)), \;\;\;\forall \;a, b, c\in A.
\end{eqnarray*}
(ii) If moreover $A$ has a usual unit $1_A$ satisfying $h\cdot 1_A=
\varepsilon (h)1_A$ for all $h\in H$, then 
\begin{eqnarray*}
&&h\cdot (ab)=(h_1\cdot a)(h_2\cdot b), 
\end{eqnarray*}
for all $h\in H$ and $a, b\in A$, that is $A$ is a left $H$-module algebra, 
so the multiplication (\ref{lulu}) is just the one of the smash product 
$A\# H$.
\end{lemma}
\begin{proof}
(i) Let $a, b, c\in A$; then one can easily compute that in $A\ot H$ we have:
\begin{eqnarray*}
&&((a\ot 1)(b\ot 1))(c\ot 1)=(y^1\cdot ((x^1\cdot a)(x^2\cdot b)))
(y^2x^3_1\cdot c)\ot y^3x^3_2,\\
&&(a\ot 1)((b\ot 1)(c\ot 1))=(y^1\cdot a)(y^2\cdot ((x^1\cdot b)(x^2\cdot c)))
\ot y^3x^3.
\end{eqnarray*}
Since $A\ot H$ is associative, these are equal; by applying $\varepsilon $ on 
the second position, we obtain 
$(ab)c=(X^1\cdot a)((X^2\cdot b)(X^3\cdot c))$, q.e.d.\\
(ii) Let $a, b\in A$ and $h\in H$; write that $((1_A\ot h)(a\ot 1))(b\ot 1)=
(1_A\ot h)((a\ot 1)(b\ot 1))$ in $A\ot H$, then apply $\varepsilon $ in the  
second position and obtain $(h_1\cdot a)(h_2\cdot b)=h\cdot (ab)$, q.e.d.
\end{proof} 
The main ingredient for proving our structure theorem for quasi-Hopf 
comodule algebras will be the structure theorem for quasi-Hopf bimodules, 
so we recall first some facts from \cite{hn}.\\
Let $H$ be a quasi-bialgebra and $M$ an $H$-bimodule together with an 
$H$-bimodule map $\rho :M\rightarrow M\ot H$, with notation 
$\rho (m)=m_{(0)}\ot m_{(1)}$ for $m\in M$ ($\rho $ is called a right 
$H$-coaction on $M$).  
Then $(M, \rho )$ is called a (right) quasi-Hopf $H$-bimodule if 
\begin{eqnarray}
&&(id_M\ot \varepsilon )\circ \rho =id_M, \label{qb1}\\
&&\Phi \cdot (\rho \ot id_M)(\rho (m))=(id_M\ot \Delta )(\rho (m))\cdot 
\Phi , \;\;\;\forall \;m\in M. \label{qb2}
\end{eqnarray}
The category of right quasi-Hopf $H$-bimodules will be denoted by 
$_H{\cal M}_H^H$ (the morphisms in the category are the $H$-bimodule maps 
intertwining the $H$-coactions).\\  
If $(V, \triangleright )$ is a left $H$-module, then $V\ot H$ becomes a   
right quasi-Hopf $H$-bimodule with structure: 
\begin{eqnarray}
&&a\cdot (v\ot h)\cdot b=(a_1\triangleright v)\ot a_2hb, \\
&&\rho _{V\ot H}(v\ot h)=(x^1\triangleright v\ot x^2h_1)\ot x^3h_2, 
\end{eqnarray}
for all $a, b, h\in H$ and $v\in V$.\\
Suppose now that $H$ is a quasi-Hopf algebra and $(M, \rho )$ is a 
right quasi-Hopf $H$-bimodule. Define the map $E:M\rightarrow M$, by 
\begin{eqnarray}
E(m)=q^1\cdot m_{(0)}\cdot \beta S(q^2m_{(1)}), \;\;\;\;\;\forall \;\;m\in M,
\end{eqnarray}
where $q=q_R=q^1\ot q^2$ is given by  
(\ref{right}). Also, for $h\in H$ and $m\in M$, define 
\begin{eqnarray}
&&h\triangleright m=E(h\cdot m). \label{act}
\end{eqnarray}
Some properties of $E$ and $\triangleright $ are collected in \cite{hn}, 
Proposition 3.4, for instance (for $h, h'\in H$ and $m\in M$): 
$E^2=E$; $E(m\cdot h)=E(m)\varepsilon (h)$; $h\triangleright E(m)=
E(h\cdot m)\equiv h\triangleright m$; $(hh')\triangleright m=h\triangleright  
(h'\triangleright m)$; $h\cdot E(m)=(h_1\triangleright E(m))\cdot h_2$; 
$E(m_{(0)})\cdot m_{(1)}=m$; $E(E(m)_{(0)})\ot E(m)_{(1)}=E(m)\ot 1$. \\
Because of these properties, the following notions of {\it coinvariants} all 
coincide: 
\begin{eqnarray*}
&&M^{co (H)}=E(M)=\{n\in M/E(n)=n\}=\{n\in M/E(n_{(0)})\ot n_{(1)}=E(n)
\ot 1\}.
\end{eqnarray*}
From the above properties it follows that $(M^{co(H)}, \triangleright )$ is 
a left $H$-module.\\
Another description of $M^{co (H)}$ is (\cite{hn}, Corollary 3.9):
\begin{eqnarray*}
&&M^{co (H)}=\{n\in M/\rho (n)=(x^1\triangleright n)\cdot x^2\ot x^3\}.
\label{co}
\end{eqnarray*}   
For a quasi-Hopf $H$-bimodule of type $V\ot H$, with $V\in $$\;_H{\cal M}$, 
we have $(V\ot H)^{co (H)}=V\ot 1$ and $E(v\ot h)=v\ot \varepsilon (h)1$,   
for all $v\in V$ and $h\in H$.\\
We can state now the structure theorem for quasi-Hopf $H$-bimodules. 
\begin{theorem} (\cite{hn}) Let $H$ be a quasi-Hopf algebra and $M$ a 
right quasi-Hopf $H$-bimodule. Consider $V=M^{co (H)}$ as a left $H$-module 
with $H$-action $\triangleright $ as in (\ref{act}), and $V\ot H$ as a right 
quasi-Hopf $H$-bimodule as above. Then the map 
\begin{eqnarray}
&&\nu :V\ot H\rightarrow M,\;\;\;\nu (v\ot h)=v\cdot h, \;\;\;\forall \; 
v\in V \;and\; h\in H, 
\end{eqnarray}
provides an isomorphism of right quasi-Hopf $H$-bimodules, with inverse 
\begin{eqnarray}
&&\nu ^{-1}:M\rightarrow V\ot H, \;\;\;\nu ^{-1}(m)=E(m_{(0)})\ot m_{(1)},  
\;\forall \;m\in M.
\end{eqnarray}
\end{theorem} 
From now on, we fix a quasi-Hopf algebra $H$ and a right $H$-comodule 
algebra $(B, \rho , \Phi _{\rho })$, with notation $\rho (b)=b_{(0)}\ot 
b_{(1)}\in B\ot H$, such that there exists $v:H\rightarrow B$ a morphism of 
right $H$-comodule algebras (in particular, this implies $\rho (v(h))=
v(h_1)\ot h_2$, for all $h\in H$, and  
$\Phi _{\rho }=v(X^1)\ot X^2\ot X^3$).   
\begin{lemma}\label{l2}
$(B, \rho )$ becomes an object in $_H{\cal M}_H^H$. 
\end{lemma} 
\begin{proof}
First, $B$ becomes an $H$-bimodule via $v$ (i.e. $h\cdot b\cdot h'=
v(h)bv(h')$ for all $h, h'\in H$ and $b\in B$). 
We prove now that $\rho :B\rightarrow B\ot H$ is an $H$-bimodule map. We 
compute:
\begin{eqnarray*}
\rho (h\cdot b\cdot h')&=&\rho (v(h)bv(h'))\\
&=&\rho (v(h))\rho (b)\rho (v(h'))\\
&=&(v(h_1)\ot h_2)(b_{(0)}\ot b_{(1)})(v(h'_1)\ot h'_2)\\
&=&v(h_1)b_{(0)}v(h'_1)\ot h_2b_{(1)}h'_2\\
&=&h_1\cdot b_{(0)}\cdot h'_1\ot h_2b_{(1)}h'_2\\
&=&h\cdot \rho (b)\cdot h', \;\;\;q.e.d.
\end{eqnarray*}
Obviously we have $(id_B\ot \varepsilon )\circ \rho =id_B$. Finally, it is  
easy to see that   
\begin{eqnarray*}
&&\Phi \cdot (\rho \ot id_B)(\rho (b))=(id_B\ot \Delta )
(\rho (b))\cdot \Phi ,  
\end{eqnarray*}
because this is exactly the condition 
\begin{eqnarray*}
&&\Phi _{\rho }(\rho \ot id_B)(\rho (b))=
(id_B\ot \Delta )(\rho (b))\Phi _{\rho } 
\end{eqnarray*}
from the definition of a right $H$-comodule algebra, due to the fact that 
$\Phi _{\rho }=v(X^1)\ot X^2\ot X^3$. Hence $(B, \rho )$ is indeed a  
right quasi-Hopf $H$-bimodule. 
\end{proof}
Since $B$ is an object in $_H{\cal M}_H^H$, we can consider the map 
$E:B\rightarrow B$, which is given by  
\begin{eqnarray}
&&E(b)=v(q^1)b_{(0)}v(\beta S(q^2b_{(1)})), \;\;\;\forall \;b\in B, 
\end{eqnarray}
where $q=q_R=q^1\ot q^2$ is given by (\ref{right}), and we can take the 
coinvariants 
\begin{eqnarray}
&&B^{co (H)}=E(B)=\{b\in B/E(b)=b\}=\{b\in B/E(b_{(0)})\ot b_{(1)}=
E(b)\ot 1\}.
\end{eqnarray}
The $H$-module algebra $A$ we are looking for will be, as a vector space, 
$A=B^{co (H)}$. \\
We have the $H$-action on $B$ given by $h\triangleright b=E(v(h)b)$, which   
gives a left $H$-module structure on $A$. Let us note that, because of 
(\ref{cucu}), we have $E(1)=1$, hence $1\in A$. 
By the structure theorem for quasi-Hopf $H$-bimodules,  
we know that the map 
\begin{eqnarray}
&&\Psi :A\ot H\rightarrow B, \;\;\;\Psi (a\ot h)=av(h), 
\end{eqnarray}
is an isomorphism in $_H{\cal M}_H^H$ (the left $H$-module  
structure of $A$ is $\triangleright $), with inverse  
\begin{eqnarray}
&&\Psi ^{-1}:B\rightarrow A\ot H, \;\;\;\Psi ^{-1}(b)=E(b_{(0)})\ot b_{(1)}.
\end{eqnarray}
Our aim will be to introduce a new multiplication on $A$, denoted by $*$,  
such that $(A, *, 1, \triangleright )$ becomes a left $H$-module algebra and 
$\Psi $ becomes an isomorphism of right $H$-comodule algebras between 
$A\# H$ and $B$ (note that $\Psi $ has the property that $\Psi \circ j=v$, 
where $j$ is the canonical map $H\rightarrow A\# H$).   
We will define actually a new multiplication $*$ on the whole $B$, and will 
take its restriction to $A$. Namely, for all $b, b'\in B$, define 
\begin{eqnarray}
&&b*b'=E(bb').
\end{eqnarray}
Since $A=E(B)$, $*$ restricts to a multiplication on $A$.  
Since for $a\in A$ we have $E(a)=a$, we obtain 
$a*1=1*a=E(a)=a$, hence $1$ is a unit for $(A, *)$. 
Let now $h\in H$; we compute:
\begin{eqnarray*}
h\triangleright 1&=&E(v(h))\\
&=&v(q^1)v(h)_{(0)}v(\beta S(q^2v(h)_{(1)}))\\
&=&v(q^1)v(h_1)v(\beta S(q^2h_2))\\
&=&v(q^1h_1\beta S(h_2)S(q^2))\\
{\rm (\ref{q5})}&=&v(q^1\beta S(q^2))\varepsilon (h)\\
{\rm (\ref{cucu})}&=&\varepsilon (h)1.
\end{eqnarray*} 
In view of Lemma \ref{l1}, in order to get that $(A, *, 1, \triangleright )$ 
is a left $H$-module algebra, it is enough to prove that the 
multiplication defined on $A\ot H$ by 
\begin{eqnarray*}
&&(a\ot h)(a'\ot h')=(x^1\triangleright a)*(x^2h_1\triangleright a')\ot  
x^3h_2h'
\end{eqnarray*}
is associative. Since $\Psi :A\ot H\rightarrow B$ is bijective and 
$B$ is associative, it is enough to prove that $\Psi $ is multiplicative, 
that is, for all $a, a'\in A$ and $h, h'\in H$:
\begin{eqnarray*}
&&\Psi ((x^1\triangleright a)*(x^2h_1\triangleright a')\ot 
x^3h_2h')=\Psi (a\ot h)\Psi (a'\ot h').
\end{eqnarray*}
We prove first a relation that will be used in the proof of the 
multiplicativity of $\Psi $.
\begin{lemma}\label{l3}
Let $H$ be a quasi-Hopf algebra; then we have: 
\begin{eqnarray}
&&q_1^1t^1x^1\ot q^1_{(2,1)}t^2_1z^1x^2\ot q^1_{(2,2)}t^2_2z^2\beta  
S(q^2t^3z^3)x^3=1\ot 1\ot 1,\label{lema3}
\end{eqnarray}
where $q=q_R=q^1\ot q^2$ is given by (\ref{right}).
\end{lemma}
\begin{proof}
We will use also the element $p=p_R=p^1\ot p^2$ given by (\ref{right}). 
We compute:\\[2mm]
${\;\;\;\;\;}$$q_1^1t^1x^1\ot q^1_{(2,1)}t^2_1z^1x^2\ot  
q^1_{(2,2)}t^2_2z^2\beta S(q^2t^3z^3)x^3$
\begin{eqnarray*}
{\rm (\ref{q3})}&=&q^1_1Z^1t^1_1y^1x^1\ot q^1_{(2,1)}Z^2t^1_2y^2x^2\ot   
q^1_{(2,2)}Z^3t^2y^3_1\beta S(q^2t^3y^3_2)x^3\\
{\rm (\ref{q5}), (\ref{q7})}&=&q^1_1Z^1t^1_1x^1\ot  
q^1_{(2,1)}Z^2t^1_2x^2\ot q^1_{(2,2)}Z^3t^2\beta S(q^2t^3)x^3\\
{\rm (\ref{q1}), (\ref{right})}&=&Z^1q^1_{(1,1)}p^1_1x^1\ot 
Z^2q^1_{(1,2)}p^1_2x^2\ot Z^3q^1_2p^2S(q^2)x^3\\
{\rm (\ref{unu})}&=&Z^1x^1\ot Z^2x^2\ot Z^3x^3\\
&=&1\ot 1\ot 1,
\end{eqnarray*}
and the relation is proved. 
\end{proof}
We will need two of the general properties of the map $E$ on a right 
quasi-Hopf $H$-bimodule $M$ recalled before, which for $M=B$ become: 
\begin{eqnarray}
&&v(h)a=(h_1\triangleright a)v(h_2),\;\;\;\forall\; 
h\in H,\;a\in A, \label{ve}
\end{eqnarray}
and, if $a\in B$, then 
\begin{eqnarray}
&&a\in A \;\Longleftrightarrow \;a_{(0)}\ot a_{(1)}=
(x^1\triangleright a)v(x^2)\ot x^3.\label{aa}
\end{eqnarray}
We need also a more explicit formula for $a*a'$, if $a, a'\in A$. 
We can write:
\begin{eqnarray*}
a*a'&=&E(aa')\\
&=&v(q^1)a_{(0)}a'_{(0)}v(\beta S(q^2a_{(1)}a'_{(1)}))\\
{\rm (\ref{aa})}&=&v(q^1)(t^1\triangleright a)v(t^2)(z^1\triangleright a')
v(z^2\beta S(q^2t^3z^3)).
\end{eqnarray*}
We can finally prove that $\Psi $ is multiplicative. We compute 
(for $a, a'\in A$ and $h, h'\in H$):
\begin{eqnarray*}
\Psi (a\ot h)\Psi (a'\ot h')&=&av(h)a'v(h')\\
{\rm (\ref{ve})}&=&a(h_1\triangleright a')v(h_2h'),
\end{eqnarray*}
${\;\;\;\;\;\;\;\;\;\;}$
$\Psi ((x^1\triangleright a)*(x^2h_1\triangleright a')\ot x^3h_2h')$
\begin{eqnarray*}
&=&(x^1\triangleright a)*(x^2h_1\triangleright a')v(x^3h_2h')\\
&=&v(q^1)(t^1x^1\triangleright a)v(t^2)(z^1x^2h_1\triangleright a')
v(z^2\beta S(q^2t^3z^3)x^3h_2h')\\
{\rm (\ref{ve})}&=&(q^1_1t^1x^1\triangleright a)v(q^1_2t^2)
(z^1x^2h_1\triangleright a')v(z^2\beta S(q^2t^3z^3)x^3h_2h')\\
{\rm (\ref{ve})}&=&(q^1_1t^1x^1\triangleright a)(q^1_{(2,1)}t^2_1z^1x^2h_1
\triangleright a')v(q^1_{(2,2)}t^2_2z^2\beta S(q^2t^3z^3)x^3h_2h')\\
{\rm (\ref{lema3})}&=&a(h_1\triangleright a')v(h_2h'),\;\;\;q.e.d. 
\end{eqnarray*}
Since obviously we have $\Psi (1\ot 1)=1$, now we have that $(A, *, 1, 
\triangleright )$ is a left $H$-module algebra and $\Psi :A\# H 
\rightarrow B$ is an algebra isomorphism. Using (\ref{aa}), the fact that 
$\rho (v(h))=v(h_1)\ot h_2$ and the formula $\rho _{A\# H}(a\# h)=
(x^1\triangleright \;a\# x^2h_1)\ot x^3h_2$, one can easily see that  
$\rho _B\circ \Psi =(\Psi \ot id)\circ \rho _{A\# H}$.  
Moreover, since 
$\Phi _{A\# H}=(1\# X^1)\ot X^2\ot X^3$ and $\Psi (1\# X^1)\ot X^2\ot X^3=
v(X^1)\ot X^2\ot X^3=\Phi _B$, we conclude that $\Psi $ is an 
isomorphism of right $H$-comodule algebras. Hence, we have proved the 
desired structure theorem:
\begin{theorem}
Let $H$ be a quasi-Hopf algebra and $B$ a right $H$-comodule algebra such 
that there exists $v:H\rightarrow B$ a morphism of right 
$H$-comodule algebras. Then there exists a left $H$-module algebra $A$ 
(whose structure is described above) such that $B\simeq A\# H$ as 
right $H$-comodule algebras.
\end{theorem}
Let now $H$, $B$ and $v:H\rightarrow B$ be as above. Since $B$ is an 
associative algebra and $v$ is an algebra map, we can consider the left 
$H$-module algebra $B^v$ as in the Preliminaries. 
\begin{proposition}
With notation as above, the map 
\begin{eqnarray*}
&&\theta :A\rightarrow B^v, \;\;\;\theta (a)=(p^1\triangleright a)v(p^2), 
\end{eqnarray*}
where $p=p_R=p^1\ot p^2$ is given by (\ref{right}), 
is an injective morphism of 
left $H$-module algebras.
\end{proposition}
\begin{proof}
Since $\Psi :A\# H\rightarrow B$ is an algebra map satisfying 
$\Psi \circ j=v$, by \cite{pvo}, Lemma 4.1. it follows that 
$\Psi :(A\# H)^j\rightarrow B^v$ is a morphism of left $H$-module algebras. 
We know from the Preliminaries that the map $i_0:A\rightarrow (A\# H)^j$, 
$i_0(a)=p^1\triangleright a\ot p^2$, is also a morphism of left $H$-module 
algebras, and one can see that actually $\theta =\Psi \circ i_0$, hence 
$\theta $ is indeed a morphism of left $H$-module algebras, and it is 
injective since $i_0$ is injective and $\Psi $ is bijective. \\
Note that in the Hopf case $\theta $ is simply the inclusion of 
$A$ into $B^v$.
\end{proof}  
\begin{remark}{\em 
Let $H$ be a quasi-Hopf algebra, $A$ a left $H$-module algebra and 
$B=A\# H$; then, for this $B$ together with the canonical map 
$j:H\rightarrow A\# H$, one can show that the map $E$ is given  
by $E(a\# h)=\varepsilon (h)(a\# 1)$, and so $B^{co (H)}=A\# 1$, with 
multiplication and $H$-action:
\begin{eqnarray*}
&&(a\# 1)* (a'\# 1)=aa'\# 1, \;\;\;\forall \;a,\;a'\in A,\\
&&h\triangleright (a\# 1)=h\cdot a\#1, \;\;\;\forall \;h\in H\;\;and
\;\;a\in A,   
\end{eqnarray*}
that is $B^{co (H)}\simeq A$ as left $H$-module algebras. Hence, the 
structure theorem allows to recover the structure of $A$ from the one 
of $A\# H$.}
\end{remark}

\end{document}